\documentclass[onefignum,onetabnum]{siamart190516}

\usepackage{amsfonts}
\usepackage{graphicx}
\usepackage{epstopdf}
\ifpdf
  \DeclareGraphicsExtensions{.eps,.pdf,.png,.jpg}
\else
  \DeclareGraphicsExtensions{.eps}
\fi

\usepackage{tikz}
\usepackage{caption}
\usepackage[labelformat=simple]{subcaption}

\usepackage{algorithmic}
\Crefname{ALC@unique}{Line}{Lines}
\usepackage{ulem}

\newsiamremark{remark}{Remark}
\newsiamremark{hypothesis}{Hypothesis}
\crefname{hypothesis}{Hypothesis}{Hypotheses}
\newsiamthm{claim}{Claim}

\newcommand{\calF}{\mathcal{F}}
\newcommand{\calL}{\mathcal{L}}
\newcommand{\calB}{\mathcal{B}}
\newcommand{\calG}{\mathcal{G}}

\newcommand{\code}[1]{\texttt{#1}}
\newcommand{\rev}[1]{{\color{blue} #1}}
\def\ExaTron{\textsc{ExaTron}}

\headers{Leveraging GPU batching for scalable nonlinear programming}{Y. Kim, F. Pacaud, K. Kim, and M. Anitescu}

\title{Leveraging GPU batching for scalable nonlinear \\ programming through massive Lagrangian decomposition\thanks{Submitted to the editors DATE.
\funding{This research was supported by the Exascale Computing Project (17-SC-20-SC),
  a collaborative effort of the U.S. Department of Energy Office of Science
  and the National Nuclear Security Administration.
  The material was based in part on work supported by the U.S. Department of Energy, Office of Science, under contract DE-AC02-06CH11357.
  This research used resources of the Oak Ridge Leadership Computing Facility at the Oak Ridge National Laboratory, which is supported by the Office of Science of the U.S. Department of Energy under Contract No. DE-AC05-00OR22725.}}}

\author{Youngdae Kim\thanks{Mathematics and Computer Science Division, Argonne National Laboratory, Lemont, IL (\email{youngdae@anl.gov}, \email{fpacaud@anl.gov}, \email{kimk@anl.gov}, \email{anitescu@mcs.anl.gov}).}
\and François Pacaud\footnotemark[2]
\and Kibaek Kim\footnotemark[2]
\and Mihai Anitescu\footnotemark[2]}

\usepackage{amsopn}

\ifpdf
\hypersetup{
  pdftitle={ExaTron},
  pdfauthor={Y. Kim, F. Pacaud, K. Kim, and M. Anitescu}
}
\fi

\begin{document}

\maketitle

\begin{abstract}
We present the implementation of a trust-region Newton algorithm \ExaTron{} for bound-constrained nonlinear programming problems, fully running on multiple GPUs.
Without data transfers between CPU and GPU, our implementation has achieved the elimination of a major performance bottleneck under a memory-bound situation, particularly when solving many small problems in batch.
We discuss the design principles and implementation details for our kernel function and core operations.
Different design choices are justified by numerical experiments.
By using the application of distributed control of alternating current optimal power flow, where a large problem is decomposed into many smaller nonlinear programs using a Lagrangian approach, we demonstrate computational performance of \ExaTron{} on the Summit supercomputer at Oak Ridge National Laboratory.
Our numerical results show the linear scaling with respect to the batch size and the number of GPUs and
more than 35 times speedup on 6 GPUs than on 40 CPUs available on a single node.
\end{abstract}

\begin{keywords}
  GPU optimization solver, second-order optimization method, decomposition
\end{keywords}

\begin{AMS}
  65K05, 90C06, 90C30, 90-04, 90-08
\end{AMS}

\section{Introduction}
\label{sec:intro}

Batch nonlinear programming refers to computing solutions of a batch of nonlinear programming problems that can be solved in parallel.
The need for  batch programming stems from distributed computing, which has been widely used to tackle large-scale optimization problems by the decomposition of the large problem into many smaller subproblems.
In this case, the subproblems are solved in parallel at each iteration of the algorithm.
Examples of distributed computing algorithms include Lagrangian-based decomposition algorithms, for example the alternating direction method of multipliers (ADMM)~\cite{Boyd11}, and value function iteration of dynamic programming~\cite{Cai15}.
Structures amenable to such decomposition are embedded in many applications in the literature including network optimization (e.g., communication networks~\cite{palomar2006tutorial}, electric grids~\cite{molzahn2017survey}, and water networks~\cite{sampathirao2017gpu}), as well as stochastic optimization~\cite{kim2018algorithmic,kim2019asynchronous} and multi-period optimization~\cite{kim2018temporal}.
In these cases, the solution time highly depends on the time to solve each batch.

While graphics processing units (GPUs) have shown  great success in accelerating the computation time of some batch operations, such as mini-batch training in machine learning and batched factorization~\cite{AbdelfattahDongarra18, HaidarDongarra18, KurzakDongarra16} in linear algebra, little attention has been given to accelerating batched nonlinear programming using GPUs.
Many studies have focused on solving a single medium to large optimization problem and leveraged GPUs to improve the computation time of the linear algebra only.
For example, GPUs have been used to accelerate the solution of linear systems arising in convex optimization algorithms~\cite{ODonoghue:2016aa,schubiger2020gpu,Smith12} and the KKT system of an augmented Lagrangian of nonlinear programming~\cite{Cao16}.
In \cite{Fei:2014aa}, the authors have implemented some number of components of the L-BFGS-B algorithm\footnote{This is a variant of the limited-memory Broyden--Fletcher--Goldfarb--Shanno (L-BFGS) for box-constrained nonlinear optimization~\cite{zhu1997algorithm}.} on GPUs, however, the algorithm does not work fully on GPUs.

All the aforementioned work takes a hybrid approach; some part requires execution on CPUs and involves data transfers between CPUs and GPUs.
Limited success in terms of outperforming CPUs has been observed on selected problems with these approaches.
This situation is partly due to lack of a fast large-scale sparse symmetric indefinite linear system solver\footnote{Currently, Nvidia's cuSOLVER~\cite{CuSolver21} library supports only a dense version of symmetric indefinite linear solver.} for GPUs, expensive data transfer cost, and the sequential nature of many optimization algorithms prohibiting full utilization of the parallel computation capability of GPUs.
Such a linear solver is required because the system of equations (the KKT matrix) to be solved at each iteration of optimization algorithms is large, sparse, symmetric, and indefinite.
In order to find a descent direction in this case, many optimization algorithms, such as Ipopt~\cite{Wachter06}, perform inertia control which requires computation of the number of positive, negative, and zero eigenvalues.
Recent experimental results~\cite{Tasseff19} of Ipopt using a GPU-based linear solver SPRAL~\cite{Duff20,Hogg16} show limited success on dense problems only.
In sparse cases, Ipopt with a GPU-based solver showed much slower performance than when it used CPU-based linear solvers.

In this paper we focus on batch nolninear programming using GPUs, where we have a large number (e.g., tens of thousands) of problems in the batch.
In this case, we employ a massive number of threads on GPUs to their full capacity to solve as many problems as possible in parallel so that we can maximize the throughput.
Although individual solve of a problem might be still slower than CPUs, a vast number of parallel computations is expected to outperform CPUs in the case of a batched solve, even CPUs with multiple cores as demonstrated in~\cref{sec:exp}.

Special attention should be given to the kernel design and its implementation according to the problem size.
When the size of problems in a batch is small, the algorithm becomes memory bound, therefore, data transfers between CPU and GPU could cause a significant cost.
Reducing such transfers is a key to accelerating computation time.
Also, due to the scarcity of available resources on the GPU hardware, kernel design should take into account factors leading to maximizing the throughput of the computation.
We discuss these factors and present our kernel design scheme in~\cref{sec:impl}.

As a GPU solver for the individual solve of a problem in a batch, we implement a novel GPU-accelerated algorithm for bound-constrained nonlinear nonconvex optimization problems of the form:
\begin{equation}
    \begin{aligned}
        & \underset{x}{\text{minimize}} && f(x) && \text{subject to} && l \le x \le u,
    \end{aligned}
    \label{eq:prob}
\end{equation}
where $x \in \mathbf{R}^d$ is the optimization variable and $l,u \in \mathbf{R}^d \cup \{-\infty,\infty\}^d$ are respectively lower and upper bounds (allowing negative and positive infinite values).
Bound constraints hold componentwise, and the objective function $f: \mathbf{R}^d \rightarrow \mathbf{R}$ is a generic nonlinear nonconvex function.
Bound-constrained problems play an important role as a building block to solve problems with more general constraints such as $h(x)=0$, where $h$ is a linear or a nonlinear function.\footnote{General inequality constraints can be formulated as equality constraints by introducing slack variables, for example, $g(x) \le 0 \Leftrightarrow g(x)+s=0, s \ge 0$.}
This is achieved by reformulating the given problem into a bound-constrained form by taking an augmented Lagrangian $L_\rho(x;\lambda,\rho) := f(x) + \lambda^T h(x) + \frac{\rho}{2}\|h(x)\|^2$ and minimizing it with respect to $x$, that is, $\min_{l \le x \le u} L_\rho(x;\lambda,\rho)$.
A sequence $\{x^k\}$ is generated with each $x^k$ corresponding to an approximate solution of the augmented Lagrangian problem, and under suitable assumptions~\cite{ConnGouldToint91, NocedalWright06} the sequence converges to a solution of the original problem with an appropriate updating scheme for $(\lambda^k,\rho^k)$.

We demonstrate the computational performance of our algorithm by solving many nonlinear programming problems of~\eqref{eq:prob} by leveraging GPU batching.
We highlight that our numerical experiment cannot be performed by any existing solver.
The solution of such nonlinear programming problems is required for the Lagrangian-based decomposition methods.
For example, in our computational experiment, the decomposition of electric grid network results in 34,704 nonlinear programming subproblems that need to be solved multiple times in the decomposition method.
Our GPU-accelerated algorithm will be used to solve the nonlinear subproblems on multiple GPUs available on the Summit supercomputer at Oak Ridge National Laboratory.
In addition, we implement the algorithm in Julia for the following reasons: portability, performance, and productivity. 
The portability is about removing the complicated process of setting compilation flags and linking proper libraries for each platform. 
A great productivity can also be achieved by Julia's high-level operations (vs. low-level operations in C) while showing as fast performance as CUDA without requiring users to optimize the code~\cite{Besard18}.

Our batch nonlinear programming solver is extremely important in the new optimization algorithm paradigm that solves a large-scale optimization by various decomposition methods (e.g.,~\cite{ryu2021privacypreserving,zhang2021tightness,SunSun21}). 
In particular, our solver enables the scalable solution of large-scale nonlinear constrained optimization problem solely on GPUs.
Recent advances on ADMM algorithms~\cite{SunSun20} enable us to have convergence guarantees even with nonconvex problems.
This implies that for any given large-scale optimization problem we can decompose it into smaller subproblems via ADMM to the extent that will work well with our batch nonlinear programming solver on GPUs.
Since other ADMM routines, such as consensus variable and multiplier updates, in general have a closed-form solution, they can be efficiently implemented on GPUs as well.
Therefore, the entire ADMM algorithm can be implemented on GPUs.
In~\cref{subsec:exp-admm}, we present such an ADMM example.

The contributions of this work can be summarized as follows:
\begin{itemize}
    \item We develop the first bound-constrained nonlinear optimization solver, implemented fully on GPU without data transfer to or from CPU.
    \item The optimization solver, as well as the ADMM used for our experiment, has been implemented in Julia, which is portable to supercomputers such as Summit at Oak Ridge National Laboratory.
    \item We have investigated and profiled multiple GPU-centric design decisions, which we report in detail in~\cref{sec:impl}. Using the optimal choices for each algorithmic unit has allowed us to obtain a superior overall performance.
    \item By applying our approach to a distributed control of large-scale electric grid operations, we demonstrate that the solution time is significantly reduced by a factor of 9--35 on 6 GPUs (vs. 40 CPU cores) on Summit.
    \item We demonstrated the multi-GPU implementation by using direct GPU-GPU communication. 
\end{itemize}

\section{Background on GPU architecture}
\label{sec:background}

\begin{figure}[htbp!]
    \centering
    \begin{subfigure}[t]{0.45\textwidth}
        \centering
        \begin{tikzpicture}
            \node at (1.7,4.2) {{\small Streaming}};
            \node at (1.7,3.9) {{\small Multiprocessor}};
            \draw (0.6,4.0) -- (0,4.0) -- (0,1.9) -- (3.5,1.9) -- (3.5,4.0) -- (2.8,4.0);

            \node at (1.7,3.3) {{\small Registers}};
            \draw (0.7,3.0) -- (2.7,3.0) -- (2.7,3.6) -- (0.7,3.6) -- (0.7,3.0);
            \draw[<->] (1.7,3.0) -- (1.7,2.7);

            \node at (1.7,2.4) {{\small L1/smem}};
            \draw (0.9,2.1) -- (2.5,2.1) -- (2.5,2.7) -- (0.9,2.7) -- (0.9,2.1);
            \draw[<->] (1.7,2.1) -- (1.7,1.7);

            \node at (1.7,1.4) {{\small L2 cache}};
            \draw (0.7,1.1) -- (2.7,1.1) -- (2.7,1.7) -- (0.7,1.7) -- (0.7,1.1);
            \draw[<->] (1.7,1.1) -- (1.7,0.8);

            \node at (1.7,0.5) {{\small DRAM}};
            \draw (0,0.2) -- (3.3,0.2) -- (3.3,0.8) -- (0,0.8) -- (0,0.2);
        \end{tikzpicture}
        \caption{Memory hierarchy}
        \label{subfig:mem}
    \end{subfigure}%
    ~
    \begin{subfigure}[t]{0.45\textwidth}
        \centering
        \begin{tikzpicture}
            \draw[->] (1.7,4.5) to [out=220,in=40] (1.7,3.8);
            \node at (2.3,4.2) {{\small Thread}};

            \node at (1.7,3.4) {{\small Thread block}};
            \draw (0.7,3.4) -- (0.5,3.4) -- (0.5,2.3) -- (2.8,2.3) -- (2.8,3.4) -- (2.7,3.4);

            \draw[->] (0.8,3.2) to [out=220,in=40] (0.8,2.5);
            \draw[->] (1.0,3.2) to [out=220,in=40] (1.0,2.5);
            \draw[->] (1.2,3.2) to [out=220,in=40] (1.2,2.5);
            \draw[->] (1.4,3.2) to [out=220,in=40] (1.4,2.5);
            \node at (1.7,2.8) {{\small ..}};
            \draw[->] (2.0,3.2) to [out=220,in=40] (2.0,2.5);
            \draw[->] (2.2,3.2) to [out=220,in=40] (2.2,2.5);
            \draw[->] (2.4,3.2) to [out=220,in=40] (2.4,2.5);

            \node at (1.7,1.9) {{\small Grid of thread blocks}};
            \draw (0.1,1.9) -- (0,1.9) -- (0,0.7) -- (3.4,0.7) -- (3.4,1.9) -- (3.3,1.9);

            \draw[->] (0.25,1.6) to [out=220,in=40] (0.25,0.9);
            \draw[->] (0.4,1.6) to [out=220,in=40] (0.4,0.9);
            \node at (0.6,1.2) {{\small ..}};
            \draw[->] (0.8,1.6) to [out=220,in=40] (0.8,0.9);
            \draw (0.1,1.65) -- (0.1,0.8) -- (1.0,0.8) -- (1.0,1.65) -- (0.1,1.65);

            \draw[->] (1.15,1.6) to [out=220,in=40] (1.15,0.9);
            \draw[->] (1.3,1.6) to [out=220,in=40] (1.3,0.9);
            \node at (1.6,1.2) {{\small ..}};
            \draw[->] (1.9,1.6) to [out=220,in=40] (1.9,0.9);
            \draw (1.0,1.65) -- (1.0,0.8) -- (2.1,0.8) -- (2.1,1.65) -- (1.0,1.65);

            \node at (2.2,1.2) {{\small ..}};

            \draw[->] (2.45,1.6) to [out=220,in=40] (2.45,0.9);
            \draw[->] (2.6,1.6) to [out=220,in=40] (2.6,0.9);
            \node at (2.9,1.2) {{\small ..}};
            \draw[->] (3.15,1.6) to [out=220,in=40] (3.15,0.9);
            \draw (2.3,1.65) -- (2.3,0.8) -- (3.3,0.8) -- (3.3,1.65) -- (2.3,1.65);

        \end{tikzpicture}
        \caption{Kernel grid}
        \label{subfig:grid}
    \end{subfigure}
    \caption{Memory hierarchy and kernel grid of GPUs.}
    \label{fig:gpu}
\end{figure}
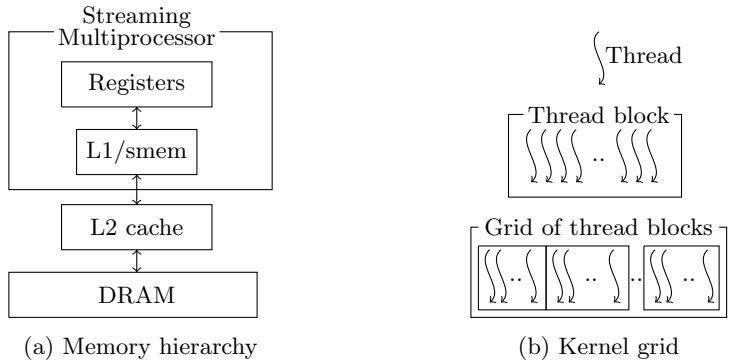

This section presents background on the GPU architecture with a focus mainly on the concepts and terminologies relevant to understanding this paper.
Note that the concepts described here are generic but use the terminologies for Nvidia's CUDA programming model~\cite{Volta17,Cuda21}, which was used for numerical experiments of this work.
For other architectures and programming models, we refer to~\cite{Amd21,Intel21}.

For GPU memory architectures, each device has its own memory hardware that is separate from the host (CPU) memory. Therefore, data transfers are needed when we want to copy data from the host to the device and vice versa.
Similar to the CPU memory system, the GPU memory is hierarchical, based on data access speed, as illustrated in~\cref{subfig:mem}.
Basically, the higher the units are in the hierarchy, the faster they are for read and write operations.
Registers are the fastest for read and write, L1 cache and shared memory are the next, and so on.
We note that the L1 cache and shared memory (smem) share the same memory unit, and that we can determine how much space we allocate for each of them.

State-of-the-art GPU programming models (such as Nvidia CUDA, OpenCL, and oneAPI) implement the execution model based on single instruction, multiple threads (SIMT) that executes single instruction on multiple threads in lockstep. 
A kernel function is a small program that is used to execute instructions on GPUs.
Once the kernel function is launched, the execution environment creates a grid of thread blocks, each of which consists of the same number of threads, as depicted in~\cref{subfig:grid}.

Each GPU architecture consists of multiple streaming multiprocessors (SMs), each of which contains scheduler and cores for computation and registers and L1 cache/shared memory for storage.
In particular, a SM schedules the execution of the instruction into warps, each of which typically consists of a set of 32 threads.
The L1 cache and shared memory are shared by all the threads inside the same thread block to help them communicate with each other, whereas registers are allocated for each thread.
Since these computational and storage resources on each SM are scarce, the number of thread blocks that can run in parallel on each SM is restricted depending on how many resources each thread block requires.
\section{Overview of the ExaTron algorithm}
\label{sec:arch}

\begin{algorithm}[t]
\begin{algorithmic}[1]
    \REQUIRE{$x^0, [l,u], \delta > 0$, and functions to evaluate $f, \nabla f$, and $\nabla^2 f$.}
    \ENSURE{$x^*$ a solution to~\eqref{eq:prob}.}
    \STATE{Set $k \leftarrow 0$ and check convergence at $x^k$.}
    \WHILE{not converged}
        \STATE\label{line3}{Evaluate $f(x^k), \nabla f(x^k)$, and $\nabla^2 f(x^k)$.}
        \STATE\label{line4}{Compute a Cauchy point: $x^k_c \leftarrow P_{[l,u]}[x^k - \alpha_k\nabla f(x^k)]$.}
        \STATE\label{line5}{Identify a subspace $\calF$ to optimize: $\calF \leftarrow \{i : l_i < (x^k_c)_i < u_i\}$}
        \STATE\label{line6}{Optimize over the subspace using a preconditioned conjugate gradient method with
        $A = \nabla^2 f(x^k)_{\calF,\calF}$:}
        \STATE\hspace{\algorithmicindent}{Compute a preconditioner $L$ such that $A + \alpha I = LL^T$ for some $\alpha \ge 0$.}
        \STATE\hspace{\algorithmicindent}{Solve $\hat{A}s = \hat{b}$ satisfying $\|s\| \le \delta$ using the conjugate gradient method where $\hat{A}=L^{-1}AL^{-T}$ and $\hat{b}= L^{-1}\nabla f(x^k)_{\calF}$.}
        \STATE\label{line9}{Set $w \leftarrow L^T\backslash s$.}
        \STATE\label{line10}{Perform a projected line search: $x^{k+1} \leftarrow P_{[l,u]}[x^k+\beta_k w]$.}
        \STATE{Set $k \leftarrow k+1$ and check convergence at $x^k$.}
    \ENDWHILE
    \RETURN{$x^* \leftarrow x^k$.}
    \caption{\ExaTron{}'s algorithm~\cite{LinMore99}}
    \label{alg:exatron}
\end{algorithmic}
\end{algorithm}

In this section we describe a variant of the trust-region Newton algorithm TRON for bound-constrained, nonlinear nonconvex optimization problems~\cite{LinMore99}.
Our algorithmic variant \ExaTron{} implements the complete Cholesky factorization for preconditioning, as opposed to the incomplete factorization in the original algorithm~\cite{LinMore99}, which is further discussed in~\cref{subsec:cholesky}.
Moreover, we present the core operations relevant to its implementation on GPUs.

\Cref{alg:exatron} summarizes the algorithmic steps of \ExaTron{}.
It is an iterative Newton-based algorithm that requires  evaluating the values of the objective function, its gradient, and Hessian at the current point at each iteration (in line 3).
These values are used to formulate a second-order Taylor approximation of the given model at the current point $x^k$, $f(x) \approx \tilde{f}(x):= f(x_k) + \nabla f(x_k)(x-x_k) + \frac{1}{2}(x-x_k)^T \nabla^2 f(x_k)(x-x_k)$.
\ExaTron{} then computes a Newton direction for the approximation $\tilde{f}$ (or a descent direction following a negative curvature in the nonconvex case) within the trust-region, expecting that we may be able to reduce the actual objective function value by moving along that direction.
These steps correspond to~\cref{line3} through~\cref{line10} of~\cref{alg:exatron}.

\ExaTron{} computes a Newton direction by the following four steps:
(i) Cauchy point computation (\cref{line4});
(ii) identification of a subspace to optimize (\cref{line5});
(iii) subspace optimization using a conjugate gradient step (\cref{line6}); and
(iv) projected line search (\cref{line10}).
The Cauchy point computation step performs the gradient projection to find a point with a sufficient reduction for the approximation $\tilde{f}$ in~\cref{line4}.
The sufficient reduction of a Cauchy point is a key property to guarantee the global convergence of the procedure.
Starting from the Cauchy point, we identify a subset of the variables to optimize further, which has been shown in~\cite{NocedalWright06,LinMore99,LinMore99-2} to result in a faster convergence rate (superlinear or even quadratic).
The subset corresponds to the variables with their values being strictly within their bounds, denoted by $\calF$ in~\cref{line5}.
We then optimize the variables in the subset by using the trust-region based conjugate gradient method~\cite{Steihaug83}.
This involves the computation of a preconditioner $L$~\cite{LinMore99-2} to accelerate the convergence of the conjugate gradient step.
Convergence of the conjugate gradient step is reached when one of the following
three conditions is satisfied:
(i) $|\calF|$ number of CG iterations has been taken;
(ii) the current CG direction reaches the trust-region radius; and
(iii) a negative curvature is detected.
We note that condition (iii) allows us to move along a descent direction even if $\nabla^2 f(x^k)$ is not positive definite~\cite{Steihaug83}.
Using the direction $w$ computed in~\cref{line9}, we perform a projected line search to find the next point with a sufficient reduction, as depicted in~\cref{line10}.
The projected line search allows a rapid change of active sets, which further accelerates the convergence.
Interested readers are referred to~\cite{LinMore99} for more details.

\Cref{fig:core-ops} describes the core operations to implement~\cref{alg:exatron}, with a mapping showing where those operations are used in it.
The numbers in parentheses at the end of each operation correspond to the line numbers of~\cref{alg:exatron}.
\ExaTron{} implements its entire algorithm on GPUs in a way that does not require data transfers between CPUs and GPUs, provided that evaluation functions are implemented on GPUs.
This removes latency issues incurred by such data transfers and enables us to achieve a much faster performance than on CPUs as described in~\cref{sec:exp}.

\begin{figure}[t]
\centering
    \begin{tikzpicture}
        \node at (1.9,3.4) {Linear algebra routines};
        \draw (0.1,3.4) -- (-0.6,3.4) -- (-0.6,0.0) -- (4.4,0.0) -- (4.4,3.4) -- (3.7,3.4);
        \node at (0.2,3.0) {{\small \code{axpy():}}};
        \node at (2.1,3.0) {{\small $y \leftarrow y + \alpha x$ (8,10)}};
        \node at (0.25,2.6) {{\small \code{ccf()}:}};
        \node at (2.1,2.6) {{\small $A +\alpha I = LL^T$ (7)}};
        \node at (0.2,2.2) {{\small \code{copy()}:}};
        \node at (1.85,2.2) {{\small $y \leftarrow x$ (4,8,10)}};
        \node at (0.27,1.8) {{\small \code{dot()}:}};
        \node at (2.1,1.8) {{\small $\sum_{i=1}^n x_iy_i$ (4,8,10)}};
        \node at (0.2,1.4) {{\small \code{gemv()}:}};
        \node at (2.45,1.4) {{\small $y \leftarrow \alpha Ax + \beta y$ (4,8,10)}};
        \node at (0.2,1.0) {{\small \code{nrm2()}:}};
        \node at (1.55,1.0) {{\small $\|x\|_2$} (4,8)};
        \node at (0.2,0.6) {{\small \code{scal()}:}};
        \node at (1.75,0.6) {{\small $x \leftarrow \alpha x$ (4,8)}};
        \node at (0.15,0.2) {{\small \code{trtrs()}:}};
        \node at (2.6,0.2) {{\small $x \leftarrow A \backslash x$ or $A^T \backslash x$ (8,9)}};

        \node at (7.0,3.4) {Other routines};
        \draw (5.8,3.4) -- (5.5,3.4) -- (5.5,0.0) -- (10.1,0.0) -- (10.1,3.4) -- (8.2,3.4);
        \node at (6.95,3.0) {{\small \code{breakpt():} (4,10)}};
        \node at (8.0,2.4) {{\small $b_i \leftarrow \left\{\begin{array}{l}
             (u_i - x_i) / w_i, w_i > 0\\
             (l_i - x_i) / w_i, w_i < 0 \end{array}\right.$}};
        \node at (7.5,1.8) {{\small \code{ccfs(): $A \leftarrow A + \alpha I$} (7)}};
        \node at (6.85,1.4) {{\small \code{gpstep():} (4,10)}};
        \node at (7.9,1.0) {{\small $s^k  \leftarrow P_{[l,u]}(x^k+\alpha_k w^k) - x^k$}};
        \node at (6.65,0.6) {{\small \code{trqsol():} (8)}};
        \node at (7.8,0.2) {{\small $\sigma$ satisfies $\|x^k + \sigma w^k\|_2=\delta$}};
    \end{tikzpicture}
    \caption{Core operations}
    \label{fig:core-ops}
\end{figure}

\section{Kernel design principles}
\label{sec:impl}

This section discuss the six kernel design principles and optimization techniques that are considered for the efficient implementation of~\cref{alg:exatron} on GPU.
The design choices made in this section are specific for a nonlinear programming problem with 32 variables or less, considering the computational experiment in~\cref{sec:exp}.
Note, however, that for problems with a larger number of variables, we discuss alternative choices of each design principle.

\subsection{Single kernel vs multiple kernels}
\label{subsec:kernel-config}

Our choice of the kernel design in \ExaTron{} is to use a single kernel function in order to maximize the algorithm throughput.
The throughput can be maximized by promptly reassigning the resources on SMs to thread blocks used for the optimization problems.
Moreover, with this design choice, a batch run can be easily implemented by launching the kernel with a grid of thread blocks, whose size is equal to the number of problems in the batch.
Assigning a specific thread block to each problem is logical, because the thread blocks are independently used and scheduled in the single-kernel design.
However, a potential caveat of this single-kernel approach is that a large number of live registers are required for a thread block, as the registers may need to save many states in different device functions during the kernel execution.

Alternatively, multiple kernels may be employed as an alternative design choice in order to avoid the high register usage issue.
However, the multi-kernel design would involve frequent memory operations to copy the current algorithm state information from one kernel to the other, exacerbating memory-bound situations.
Moreover, the implementation becomes more complicated with multiple kernel functions, because we need to track down each problem's state outside of the kernel function.
If the same kernel function was used for different problems by employing multiple thread blocks, the throughput would be significantly degraded due to the potential wait time for reassigning problems to the kernel function.

\subsection{Thread configuration}
\label{subsec:thread-config}

\begin{figure}[t]
    \centering
    \includegraphics[scale=.6]{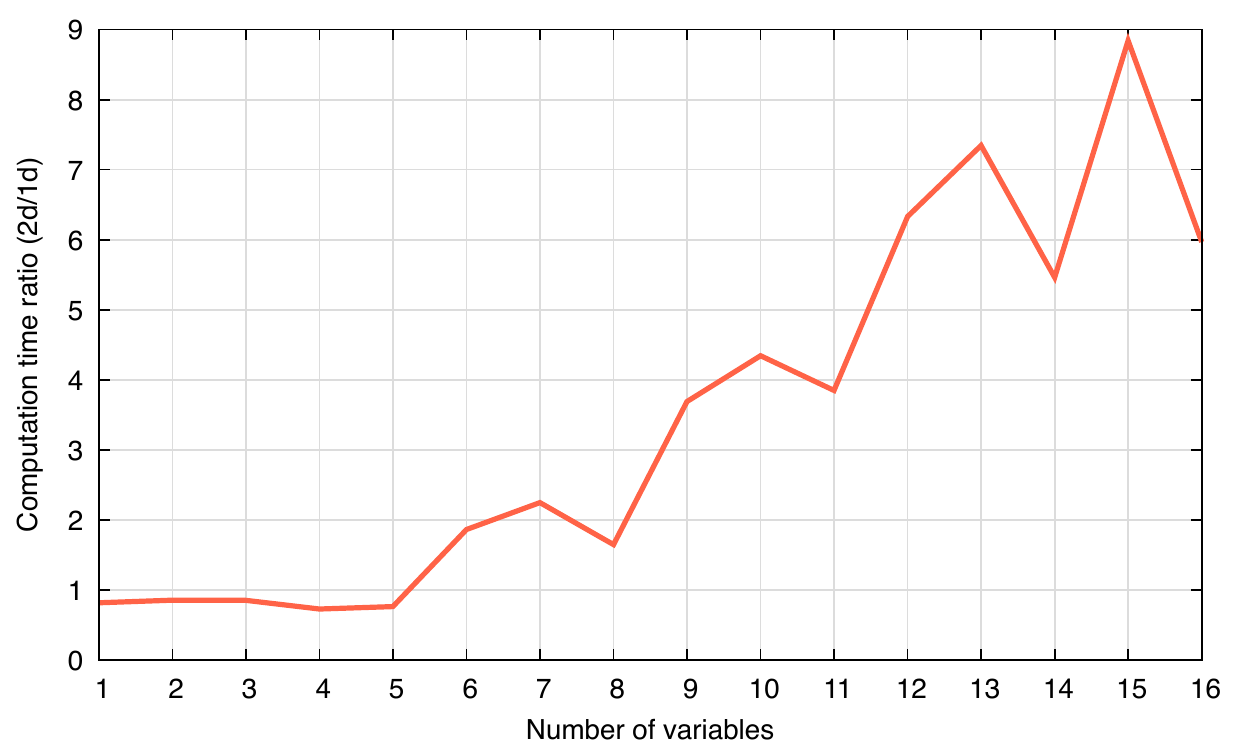}
    \caption{Computation time ratio between the 1-dimensional and 2-dimensional
             thread configurations over a batch solve of 10,000 nonlinear
             nonconvex problems.}
    \label{fig:thread-conf}
\end{figure}

\ExaTron{} uses the single warp 1-dimensional thread configuration.
Each thread block consists of a single warp of 32 threads with an 1-dimensional thread configuration.
These threads naturally enable parallel operations on a vector of size up to 32.
Since most of our operations are matrix-vector and vector-vector operations (described later in this section), such a 1-dimensional scheme fits well with our needs.

Since we perform factorization of a matrix of size $n \times n$ in~\cref{alg:exatron}, one might think that a 2-dimensional $n \times n$ thread configuration may be more efficient, with each thread in charge of each element of a matrix.
This configuration may be effective when we read or write the entire matrix: one line of code will do the work, such as $B[tx,ty] \leftarrow A[tx,ty]$, where $tx$ and $ty$ are thread IDs in $x$- and $y$-coordinate, respectively.

However, the number of threads required for the 2-dimensional $n \times n$ configuration is proportional to the square of the number of variables, limiting the throughput of \ExaTron{} on each multiprocessor significantly.
As discussed in~\cite{AbdelfattahDongarra18}, for $n=16$ we already need 256 threads for each thread block in this case.
On the Nvidia's Volta architecture~\cite{Volta17}, this implies that we can have at most 8 active thread blocks on each multiprocessor, since a maximum of 2,048 threads is allowed per multiprocessor.
Since a maximum of 32 thread blocks is allowed on each multiprocessor, the maximum throughput can be reduced by a factor of 4 with this configuration.

Another issue with the 2-dimensional $n \times n$ configuration is that instruction efficiency is expected to be much lower than the single warp 1-dimensional configuration.
Among the core operations of \ExaTron{} (\cref{fig:core-ops}), factorization and triangular solves are the most expensive.
However, they factorize or solve a single column of size at most $n$ at a time, which makes the remaining threads idle.
Moreover, the synchronization burden increases because a thread-level synchronization should be performed to synchronize between warps instead of a lighter warp-level synchronization.

In~\cref{fig:thread-conf}, we plot the computation time ratio between the 1-dimensional and 2-dimensional thread configurations obtained from running \ExaTron{} for a batch solve of 10,000 nonlinear nonconvex problems.\footnote{We modified the \texttt{hs45} problem from~\cite{Gould03} to experiment with different numbers of variables. The problem is a bound-constrained nonlinear nonconvex problem with objective function $f(x)=120 - \prod_{i=1}^n x_i$ and bounds $x_i \le i$ for $i=1,\dots,n$.
Hence, the optimal solution is $x^*_i = i$. In each batch, we generated the same 10,000
instances.}
When a single warp was used in both cases for $n \le 5$, the 2-dimensional configuration showed slightly faster computation time than did the 1-dimensional configuration.
For $n \ge 6$, however, its computation time deteriorated sharply when it started to use more than a single warp; and it showed up to 9 times slower performance than the 1-dimensional configuration.
A similar result was obtained in~\cite{AbdelfattahDongarra18} for Cholesky, LU, and QR factorizations under the same two different thread configurations.

\subsection{Cholesky factorization}
\label{subsec:cholesky}

We use the Cholesky factorization to compute a preconditioner $L$ for a dense Hessian matrix $A=\nabla^2 f(x)$.\footnote{There are other methods to compute a preconditioner, such as LU factorization. However, Cholesky factorization exploits the symmetricity of the Hessian matrix, enabling two times faster computation time than LU factorization.}
If $A$ is positive definite, the Cholesky factorization will give us a lower triangular factor $L$ such that $A=LL^T$.
In this case, $\hat{A}=L^{-1}AL^{-T}$ becomes an identity matrix, and we can find a solution to $\hat{A}s=\hat{b}$ in just one conjugate gradient iteration.\footnote{We note that because of the numerical limitation of finite precision of floating-point numbers we determine $s$ is a solution to $\hat{A}s=\hat{b}$ if $\|\hat{A}s - \hat{b}\| \le \epsilon$ for some small error tolerance $\epsilon > 0$.}
If $A$ is not positive definite, $L$ is computed for a diagonally perturbed matrix $(A + \alpha I)$ for some $\alpha > 0$.
Such a  preconditioner is expected to have eigenvalues of $\hat{A}$ clustered so that we can find a solution in a few conjugate gradient iterations~\rev{\cite{LinMore99-2,LinMore99,NocedalWright06}}.

In contrast to the existing method~\cite{LinMore99}, which implements an incomplete Cholesky factorization with sparse linear algebra, \ExaTron{} implements a complete Cholesky factorization with dense linear algebra.
The incomplete Cholesky factorization stores only a subset of the new nonzero entries generated during factorization because of storage limitation with sparse linear algebra.
This may result in a less accurate preconditioner.
In contrast, the complete Cholesky factorization stores all of the newly generated nonzero entries.
Since we are dealing with small sparse Hessian matrices, storing them in a dense matrix does not induce any memory limitations on GPUs.
Dense linear algebra computations combined with GPU's SIMT capability may be as competitive as sparse linear algebra in this case as well.
By storing all the newly generated nonzero entries, we will have a more accurate preconditioner than that of the incomplete Cholesky factorization, which could potentially leads to fewer conjugate gradient iterations.

By taking advantage of a small matrix, we store its entire elements in shared memory during factorization for efficient read/write operations on GPUs.
In this case, we do not have to employ a blocking algorithm that stores a block of elements at a time and factorizes it for efficient data reuse.
All the elements are already available in shared memory, and read/write operations on shared memory are much faster than on global memory.
Therefore, we implement an unblocked version of Cholesky factorization;
we factorize one column at a time without employing level 3 BLAS operations.
A similar design choice was made in~\cite{AbdelfattahDongarra18} for a batch of Cholesky factorizations of small matrices of size $n \le 32$.

To further optimize the implementation, we have experimented with two representative Cholesky factorization algorithms~\cite{HaidarDongarra18,KurzakDongarra16}:
left-looking factorization and right-looking factorization.
The left-looking factorization applies all the previous updates just before a column is factorized, called lazy update, whereas the right-looking factorization immediately applies the updates to the trailing submatrix right after a column is factorized.
A theoretical analysis in~\cite{HaidarDongarra18} shows that the left-looking algorithm has a smaller number of read/write operations than that of the right-looking algorithm, thus making it more efficient under memory-bound situations.

\Cref{tbl:cholesky} shows the average computation time of the left-looking and right-looking batched
Cholesky factorizations over 10,000 randomly generated matrices for $1 \le n \le 32$.
As expected from the analysis, the left-looking algorithm showed a faster performance than that of the right-looking algorithm.
Hence, \ExaTron{} uses by default the left-looking Cholesky factorization algorithm to compute a preconditioner.

\begin{table}[h]
    \caption{Average computation time of the left- and right-looking Cholesky factorizations.}
    \label{tbl:cholesky}
    \centering
    \begin{tabular}{c|c}
        Factorization & Time (ms)\\\hline
        Left-looking & 0.362\\
        Right-looking & 0.437
    \end{tabular}
\end{table}

\subsection{Triangular solves}
\label{subsec:triangular-solve}

In the conjugate gradient step, we need to perform matrix-vector multiplications between a preconditioned matrix $\hat{A}=L^{-1}AL^{-T}$ and a vector $p$.
This involves two triangular solves:
1) $L^Tz=p$ ({\it backward substitution}); and
2) $Lq=\tilde{z}$ ({\it forward substitution}), where $\tilde{z}=Az$.

Triangular solves are inherently memory bound because they have to sequentially access one column at a time with a simple arithmetic applied to each element of it.
Moreover, the number of elements to access is diminishing as we move forward or backward in both forward and backward substitutions.
Hence, at most $n$ threads are needed.
In this case, we apply data parallelism as described in~\cref{subsec:dlp}, updating each element of a column in parallel.

Implementing the forward substitution is  straightforward; however, the backward substitution involves the transpose of $L$, requiring a row-wise memory access with as many strides as the size of the matrix.
Although we store the elements of $L$ in shared memory, this memory access pattern could cause bank conflicts.
For example, when $n=8$ and the shared memory has 32 banks, a 2-way bank conflict will occur: elements in the same row with column indices differing by 4 will be stored in the same bank as they are accessed at the same time by
different threads.

To avoid bank conflicts, we store the elements for both $L$ and $L^T$ in the same matrix.
Since $L$ is a lower triangular matrix stored in a dense matrix format, we could use its upper triangular part to store the elements of $L^T$.
This would involve additional write operations for the subspace optimization step of~\cref{alg:exatron}.
But since the backward substitution is applied multiple times during the conjugate gradient iterations and we could avoid bank conflicts, we may obtain performance gain offsetting the cost of additional writes.
A similar approach was applied in~\cite{HaidarDongarra18}.

In~\cref{tbl:trsolve} we present the average computation time of explicit and implicit backward substitutions over 10,000 randomly generated matrices.
By explicit, we mean we save the transpose of $L$ in its upper triangular part explicitly, whereas we perform row-wise access in the case of implicit.
The results demonstrate that saving the transpose $L^\top$ explicitly yields a faster computation time.

\begin{table}[h]
    \caption{Average computation time of explicit and implicit triangular solve}
    \label{tbl:trsolve}
    \centering
    \begin{tabular}{c|c}
        Backward substitution & Time (ms)\\\hline
        Explicit & 0.875\\
        Implicit & 1.408
    \end{tabular}
\end{table}

\subsection{Shared-memory management}
\label{subsec:smem}

In many parts of \ExaTron{}, memory space is needed to share data between device functions and to store and reuse intermediate computation results.
The current iterate $x^k$ is used by all components of \ExaTron{} together with its lower and upper bounds.
A preconditioner matrix is referred to in multiple places of \ExaTron{}.
We also need some number of vectors to store a Cauchy point and perform the conjugate gradient step.

One advantage of \ExaTron{}'s algorithm is that it does not require dynamic memory allocation during its procedure. Therefore,  we can predetermine the memory space needed for its computation for a given problem size $n$.
For efficient data access and reuse, we use shared memory for such memory space.

\Cref{tbl:smem} presents the effect of shared memory on the computation time of \ExaTron{} over the batch of 10,000 nonlinear nonconvex problems used in~\cref{subsec:thread-config}.
With shared memory, \ExaTron{}'s computation performance was about  twice faster.\footnote{As noted on page 17 of~\cite{Volta17}, the Volta architecture significantly improved L1 data cache performance, leading to much lower latency and higher bandwidth. When we do not use shared memory, frequently used data are cached in the L1 data cache.}

\begin{table}[h]
    \caption{Average computation time with/without using shared memory}
    \label{tbl:smem}
    \centering
    \begin{tabular}{c|r}
        Data share       & Time (ms) \\ \hline
        Shared memory    & 9.536     \\
        No shared memory & 16.293
    \end{tabular}
\end{table}

Currently, the memory space allocated in shared memory for each thread block is proportional to the square of the problem size, $O(n^2)$.
Since 96~KB of shared memory are available on recent GPU architectures, for small $n$ our shared-memory requirement is not a limiting factor for achieving high occupancy.
For medium or larger size $n > 32$, however, it could become a limiting factor, and a different implementation may be needed to achieve higher occupancy in that case.

\subsection{Register file management}
\label{subsec:reg}

\ExaTron{}'s kernel invokes multiple device functions, and some of them are called in a nested way.
This could increase the number of live registers significantly by accumulating the function calls in the nested call graph.
When we compile \ExaTron{}, the output of ``\code{ptxas -v}'' shows that our kernel function uses more than 100 registers.
If we add another kernel on top of it (such as a kernel implementing an augmented Lagrangian algorithm), the number of registers could increase even more.

The use of a large number of registers potentially hinders achieving high occupancy.
On the Volta architecture, the number of registers allowed per multiprocessor is 65,536.
Assuming that each thread is using about 100 registers, the limit on the number of active thread blocks becomes $20$.\footnote{The number was obtained from $\text{floor}(65536/(32*100),4) = 20$ where we round down the result of the division to the multiple of 4, 4 being the warp allocation granularity.}
Therefore, the theoretical occupancy we could achieve is at most $31\%$.

To achieve higher occupancy, we need to reduce the number of registers down to where it stops improving performance.
In general, there is a trade-off in restricting the number of registers.
It may incur register spills to local memory, potentially slowing performance on data access.
On the other hand, it could increase occupancy so that larger number of thread blocks could be executed simultaneously.
We control the number of registers to use by specifying explicitly the compile option \code{maxregs}.

\section{Other Implementation Details}

We describe other implementation details to optimize the efficiency of the core operations on GPU.
The implementation details considered in this section are independent to the problem size.

\subsection{Register shuffling and warp-level synchronization}
\label{subsec:warp}

In some parts of \ExaTron{}, we need to compute an aggregation over a vector, such as a sum or a maximum of the vector elements, and broadcast the result back to threads.
For example, a two-norm or a maximum-norm of a vector needs to be computed to check the violation of a trust region, $\|x\|_2 \le \delta$, or the convergence of the algorithm (measured by the maximum element of the projected gradient vector $\|\nabla_{[l,u]}f(x)\|_\infty \le \epsilon$).
Once we compute such a value, we need to broadcast it back to threads so that each thread can proceed to the next step.

\begin{figure}[htbp!]
    \centering
    \begin{tikzpicture}
        \node[circle,draw=black,fill=white] (a1) at (0.5,4.0) {};
        \node[circle,draw=black,fill=white] (a2) at (1.3,4.0) {};
        \node[circle,draw=black,fill=white] (a3) at (2.1,4.0) {};
        \node[circle,draw=black,fill=white] (a4) at (2.9,4.0) {};
        \node[circle,draw=black,fill=white] (a5) at (3.7,4.0) {};
        \node[circle,draw=black,fill=white] (a6) at (4.5,4.0) {};
        \node[circle,draw=black,fill=white] (a7) at (5.3,4.0) {};
        \node[circle,draw=black,fill=white] (a8) at (6.1,4.0) {};

        \node[circle,draw=black,fill=white] (b1) at (0.5,3.0) {};
        \node[circle,draw=black,fill=white] (b2) at (1.3,3.0) {};
        \node[circle,draw=black,fill=white] (b3) at (2.1,3.0) {};
        \node[circle,draw=black,fill=white] (b4) at (2.9,3.0) {};
        \node[circle,draw=black,fill=white] (b5) at (3.7,3.0) {};
        \node[circle,draw=black,fill=white] (b6) at (4.5,3.0) {};
        \node[circle,draw=black,fill=white] (b7) at (5.3,3.0) {};
        \node[circle,draw=black,fill=white] (b8) at (6.1,3.0) {};

        \draw[->] (a1) -- (b1); \draw[->] (a5) -- (b1);
        \draw[->] (a2) -- (b2); \draw[->] (a6) -- (b2);
        \draw[->] (a3) -- (b3); \draw[->] (a7) -- (b3);
        \draw[->] (a4) -- (b4); \draw[->] (a8) -- (b4);
        \node at (0.7,3.2) {\small \code{op}};
        \node at (1.5,3.2) {\small \code{op}};
        \node at (2.3,3.2) {\small \code{op}};
        \node at (3.1,3.2) {\small \code{op}};
        \node at (6.3,3.5) {\small \code{CUDA.shfl\_down\_sync()}};

        \node[circle,draw=black,fill=white] (c1) at (0.5,2.0) {};
        \node[circle,draw=black,fill=white] (c2) at (1.3,2.0) {};
        \node[circle,draw=black,fill=white] (c3) at (2.1,2.0) {};
        \node[circle,draw=black,fill=white] (c4) at (2.9,2.0) {};
        \node[circle,draw=black,fill=white] (c5) at (3.7,2.0) {};
        \node[circle,draw=black,fill=white] (c6) at (4.5,2.0) {};
        \node[circle,draw=black,fill=white] (c7) at (5.3,2.0) {};
        \node[circle,draw=black,fill=white] (c8) at (6.1,2.0) {};

        \draw[->] (b1) -- (c1); \draw[->] (b3) -- (c1);
        \draw[->] (b2) -- (c2); \draw[->] (b4) -- (c2);
        \node at (0.7,2.2) {\small \code{op}};
        \node at (1.5,2.2) {\small \code{op}};
        \node at (6.3,2.5) {\small \code{CUDA.shfl\_down\_sync()}};

        \node[circle,draw=black,fill=white] (d1) at (0.5,1.0) {};
        \node[circle,draw=black,fill=white] (d2) at (1.3,1.0) {};
        \node[circle,draw=black,fill=white] (d3) at (2.1,1.0) {};
        \node[circle,draw=black,fill=white] (d4) at (2.9,1.0) {};
        \node[circle,draw=black,fill=white] (d5) at (3.7,1.0) {};
        \node[circle,draw=black,fill=white] (d6) at (4.5,1.0) {};
        \node[circle,draw=black,fill=white] (d7) at (5.3,1.0) {};
        \node[circle,draw=black,fill=white] (d8) at (6.1,1.0) {};

        \draw[->] (c1) -- (d1); \draw[->] (c2) -- (d1);
        \node at (0.7,1.2) {\small \code{op}};
        \node at (6.3,1.5) {\small \code{CUDA.shfl\_down\_sync()}};

        \draw[->] (d1.east) -- (d2.west);
        \draw[->] (d1.east) to [bend right=25] (d3.west);
        \draw[->] (d1.east) to [bend right=25] (d4.west);
        \draw[->] (d1.east) to [bend right=25] (d5.west);
        \draw[->] (d1.east) to [bend right=25] (d6.west);
        \draw[->] (d1.east) to [bend right=25] (d7.west);
        \draw[->] (d1.east) to [bend right=25] (d8.west);
        \node at (6.65,0.5) {\small \code{CUDA.shfl\_sync()}};
    \end{tikzpicture}
    \caption{Sharing values between threads in the same warp via register shuffling.
    \code{op} is a binary operator such as \code{+} or \code{max}.}
    \label{fig:register-shuffle}
\end{figure}
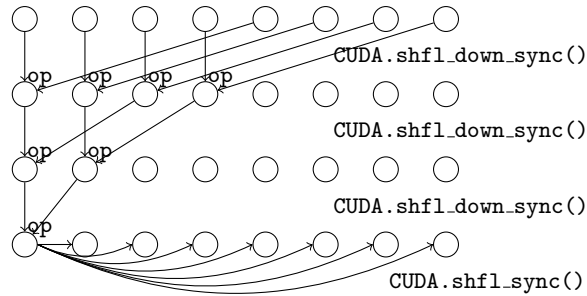

For an efficient utilization of threads and a lightweight synchronization, we implement this type of aggregation using a register shuffling and a warp-level synchronization.
The register shuffling refers to shared register values between threads in the same warp-through shuffling.
For example, in order to compute $\|x\|_2=\sqrt{\sum_{i=1}^n x_i^2}$, each thread $i$ for $i \in \{1,\dots,n\}$ computes $x_i^2$ first, and the values are summed up by using a register shuffling as shown in~\cref{fig:register-shuffle} for $n=8$.
The aggregated value is stored in the first thread, and we broadcast it back to all threads in the same warp using the register shuffling again.
Similarly, computing a maximum is performed by replacing the \code{+} operator with \code{max}.
Register shuffling routines include a warp-level synchronization such as \code{CUDA.sync\_warp()}, which is lighter than a thread-level synchronization \code{CUDA.sync\_threads()}.
Also, the use of registers avoids the use of more expensive shared memory.
All of these make it efficient to compute an aggregation and share data between threads in the same warp.

\subsection{Data-level parallelism}
\label{subsec:dlp}

Data-level parallelism (DLP) refers to applying the same instruction to multiple data in parallel.
Arithmetic operations such as \code{axpy()}, \code{copy()}, and \code{scal()} can benefit from DLP by affecting one thread per element in a vector, replacing the use of a loop by a single line of code.
For example, \code{axpy()} computes \code{y $\leftarrow$ y + ax}, where \code{x}
and \code{y} are vectors of the same size and \code{a} is a scalar.
Without the use of DLP, we have to loop around each element of those two vectors like \code{for i=1:n y[i] = y[i] + a*x[i]}, where only a single thread works at a time.
With DLP, this can be implemented as \code{y[tid] = y[tid] + a*x[tid]}, where \code{tid} is a thread ID less than or equal to $n$, utilizing $n$ threads simultaneously.
In a similar way, \ExaTron{} makes use of DLP for other arithmetic operations whenever applicable.

\subsection{Instruction-level parallelism}
\label{subsec:ilp}

Within a thread, instruction-level parallelism (ILP) refers to the simultaneous execution of multiple instructions in a sequence of them that can be executed in any order.
This depends on the compiler's capability to detect such independent instructions.

An example of ILP in \ExaTron{} can occur when we form a submatrix in step ii of the procedure.
We extract and copy the rows and the columns of a Hessian matrix $A$ corresponding to the free variables.
This is implemented by assigning each thread to each free variable index---applying DLP---and
having it copy the rows of $A$ corresponding to that column.
There are multiple rows to copy for each thread, and these copy instructions are independent of each other, making them a target to apply ILP:
for thread $j$ taking care of the $j$th free column, it copies elements $A[i,j]$ for $i \in R(j)$, with $R(j)$ being the set of row indices to copy for column $j$, and these copy instructions can be executed in any order.

\section{Computational Results}
\label{sec:exp}

In this section we demonstrate the performance of \ExaTron{} using a distributed control application from the power system literature.
In~\cref{subsec:exp-admm} we briefly introduce an alternating current optimal power flow (ACOPF) and its distributed-computation method based on an alternating direction method of multipliers (ADMM).
\Cref{subsec:exp-single-gpu,subsec:exp-multiple-gpu} describe the performance of \ExaTron{} over single and multiple GPUs, respectively.
We compare the performance of the GPU implementation of \ExaTron{} with a CPU implementation in~\cref{subsec:exp-cpu}.

The implementation of \ExaTron{} is written in Julia, and the code is available at~\url{https://github.com/exanauts/ExaTron.jl}.
All experiments were performed on a compute node of the Summit supercomputer at Oak Ridge National Laboratory using Julia 1.6.0~\cite{Julia17} and CUDA.jl 2.6.1~\cite{Besard18}.
Note, however, that our implementation is not limited to a single node.
Each compute node of the Summit supercomputer has 2 sockets of POWER9 processors having 22 physical cores each, 512 GB of DRAM, and 6 NVIDIA Tesla V100 GPUs evenly distributed to each socket~\cite{Summit21}.
We note that the MPI communication between GPUs was implemented using a CUDA-aware MPI with the asynchronous send/receive routines. 
For CPUs, we used the scatter/gather routines. Our different use of the MPI routines is mainly because we found scatter/gather routines between GPUs were much slower than the asynchronous send/receive.
We re-emphasize that \ExaTron{} does not require transferring data between the CPU and GPU.

\subsection{Distributed control of ACOPF}
\label{subsec:exp-admm}

In electrical engineering, the optimal power flow problem~\cite{WoodWollenbergSheble13} focuses on computing the optimal dispatch of active and reactive powers among a set of generators, while satisfying physical constraints such as voltage/generator limits, power balance, and transmission line limits.
The power network is modeled as a undirected graph: 
generators are located on certain buses of the network.
The active and reactive powers produced by the generators flow between buses via the transmission lines, so as to satisfy demand at each bus.
The problem corresponds to the ACOPF problem when the power flow is modeled by using alternating
current. 
In that case, there exists an exact formulation that encapsulates the actual physical constraints, such as Ohm's and Kirchhoff's current laws.
Computing an exact and efficient solution to the ACOPF problem has a  practical application, since a small percentage increase in market efficiency leads to billions of dollars of cost savings per year~\cite{Cain13} for transmission operators.

Four representative variants of the ACOPF formulations exist: polar, rectangular, extended rectangular, and current-voltage formulations.
Although their forms are different, they are equivalent between each other.
In the following, we present the extended rectangular formulation using the notation in~\cite{Mhanna19}:

\begin{subequations}
    \begin{align}
        &\underset{\substack{p_{gi},q_{gi},p_{ij},q_{ij},\\p_{ji},q_{ji},w_i,w_{ij}^R,w_{ij}^I}}{\text{minimize}} \quad \sum_{(g,i)\in\calG} f_{gi}(p_{gi}) \notag \\
        \text{sub} & \text{ject to} \notag \\
        & \underline{p}_{gi} \leq p_{gi} \leq \bar{p}_{gi}, && (g,i)\in\calG \label{eq:GenLimitReal}\\
        & \underline{q}_{gi} \leq q_{gi} \leq \bar{q}_{gi}, && (g,i)\in\calG \label{eq:GenLimitReactive}\\
        & \sqrt{p_{ij}^2 + q_{ij}^2} \leq \bar{s}_{ij},   && (i,j) \in\calL\cup\calL_t \label{eq:LineLimit}\\
        & (w_{ij}^R)^2 + (w_{ij}^I)^2 = w_i w_j, && (i,j)\in\calL \label{eq:VoltageW} \\
        & \underline{v}_i^2 \leq w_i \leq \bar{v}_i^2, && i\in\calB \label{eq:VoltageLimitRect}\\
        & \tan(\underline{\theta}_{ij}^\Delta) w_{ij}^R \leq w_{ij}^I \leq \tan(\bar{\theta}_{ij}^\Delta) w_{ij}^R, && (i,j)\in\calL \label{eq:AngleLimitRect}\\
        & \sum_{(g,i)\in\calG} p_{gi} - p_i^d = \sum_{j \in \calB_i} p_{ij} + g_i^{sh} w_i, && i\in\calB \label{eq:PowerBalanceRealRect}\\
        & \sum_{(g,i)\in\calG} q_{gi} - q_i^d = \sum_{j\in\calB_i} q_{ij} - b_i^{sh} w_i, && i\in\calB \label{eq:PowerBalanceReactiveRect}\\
        & p_{ij} = g_{ij}^c w_i - g_{ij} w_{ij}^R + b_{ij} w_{ij}^I, && (i,j) \in \calL \label{eq:pij} \\
        & q_{ij} = b_{ij}^c w_i - b_{ij} w_{ij}^R - g_{ij} w_{ij}^I, && (i,j) \in \calL \label{eq:qij} \\
        & p_{ji} = g_{ji}^c w_j - g_{ji} w_{ij}^R - b_{ji} w_{ij}^I, && (i,j) \in \calL \label{eq:pji} \\
        & q_{ji} = b_{ji}^c w_j - b_{ji} w_{ij}^R + g_{ji} w_{ij}^I, && (i,j) \in \calL, \label{eq:qji} 
    \end{align}
    \label{eq:ac-rectangular}
\end{subequations}
where $\calG, \calL$, and $\calB$ denote the set of generators, branches, and buses in the network, respectively.
We use $(g,i) \in \calG$ to represent that generator $g$ is connected to bus $i$.
For a branch $(i,j) \in \calL$, $i$ is the ``from" bus and $j$ is the ``to" bus in $(j,i) \in \calL_t$.
Formulation~\eqref{eq:ac-rectangular} is a {\it nonlinear nonconvex} optimization problem, which is known to be computationally challenging;
even verifying local optimality can be an NP-hard problem~\cite{Pardalos88}.
Hence, convergence is in general achieved at a point satisfying only second-order necessary optimality conditions~\cite{NocedalWright06}, and we are not interested in proving global optimality.

In order to efficiently solve large-scale ACOPFs, Mhanna et al.~\cite{Mhanna19} introduce a distributed control approach, where the problem is decomposed into components---generators, branches, and buses ---by duplicating the variables linking different components.
Then, the problem is solved by optimizing each component separately using an ADMM algorithm.
One subproblem is associated with each component, resulting in a total number of subproblems equal to $(|\calG|+|\calL|+|\calB|)$.
At each ADMM iteration, the algorithm starts by solving the generator and branch subproblems in parallel. 
Then, the bus subproblems are solved concurrently. 
Once all the subproblems are solved, the algorithm updates its Lagrange multipliers and moves to the next iteration.

One advantage of the algorithm introduced in~\cite{Mhanna19} is that both the generator and the bus subproblems have a closed-form solution, so we do not have to employ a nonlinear optimization solver for them.
However, the solutions of branch subproblems require to solve nonlinear nonconvex problems, each of which is formulated as~\eqref{eq:exp-branch} for branch $(i,j) \in \calL$:
\begin{subequations}
    \begin{align}
        &\underset{v_i,v_j,\theta_i,\theta_j}{\text{minimize}} &&
        \sum_{(l,m) \in \{(i,j)\cup(j,i)\}}\Bigg(
            \lambda_{p_{lm}}(p_{lm}-\tilde{p}_{lm}) +
            \lambda_{q_{lm}}(q_{lm}-\tilde{q}_{lm})
        \notag\\
        &&& + \frac{\rho_{p_{lm}}}{2}(p_{lm}-\tilde{p}_{lm})^2 +
            \frac{\rho_{q_{lm}}}{2}(q_{lm}-\tilde{q}_{lm})^2
        \Bigg) \notag\\
        &&& + \sum_{l \in \{i,j\}}\Bigg(
            \lambda_{w_{l}}(v_l^2 - \tilde{w}_l) +
            \lambda_{\theta_{l}}(\theta_l - \tilde{\theta}_l)
        \notag\\
        &&& \quad\quad\quad\quad\quad
         + \frac{\rho_{w_l}}{2}(v_l^2 - \tilde{w}_l)^2 +
         \frac{\rho_{\theta_l}}{2}(\theta_l - \tilde{\theta}_l)^2
        \Bigg)\notag\\
        &\text{subject to} &&\notag\\
        &&&\underline{v}_i \le v_i \le \overline{v}_i\\
        &&&\underline{v}_j \le v_j \le \overline{v}_j\\
        &&&-2\pi \le \theta_i, \theta_j \le 2\pi  ,
    \end{align}
    \label{eq:exp-branch}
\end{subequations}
where the definition of $p_{ij},q_{ij},p{ji}$, and $q_{ji}$ follows from \eqref{eq:pij}--\eqref{eq:qji} by plugging-in $w_i=v_i^2,w_j=v_j^2,w_{ij}^R=v_iv_j\cos(\theta_i-\theta_j)$, and $w_{ij}^I=v_iv_j\sin(\theta_i-\theta_j)$.
Note that the objective of~\eqref{eq:exp-branch} is nonconvex.

\begin{table}[htpb!]
    \caption{Data statistics}
    \label{tbl:exp-data}
    \centering
    \begin{tabular}{r|r|r|r}
        \multicolumn{1}{c|}{Data} & \# Generators & \# Branches & \# Buses\\\hline
        2868rte & 600 & 3,808 & 2868 \\
        6515rte & 1,389 & 9,037 & 6515 \\
        9241pegase & 1,445 & 16,049 & 9,241 \\
        13659pegase & 4,092 & 20,467 & 13,659\\
        19402goc & 971 & 34,704 & 19,402
    \end{tabular}
\end{table}

We have implemented the ADMM algorithm fully on GPUs without data transfer to the CPU, and we use \ExaTron{} to solve the branch subproblems at each ADMM iteration.
The ADMM algorithm has also been written in Julia.
We have experimented with our implementation in five large-scale examples from the MATPOWER~\cite{Zimmerman11} and PGLIB benchmark~\cite{babaeinejadsarookolaee2019power}, where the first four of them were also used in the literature to test ADMM algorithm~\cite{Mhanna19,SunSun21}.
\Cref{tbl:exp-data} presents the data statistics of our test examples.
We note that up to 34K nonlinear nonconvex problems are solved by \ExaTron{} at each ADMM iteration.

\subsection{Performance on a single GPU}
\label{subsec:exp-single-gpu}

\begin{figure}[htbp!]
    \centering
    \includegraphics[scale=.65]{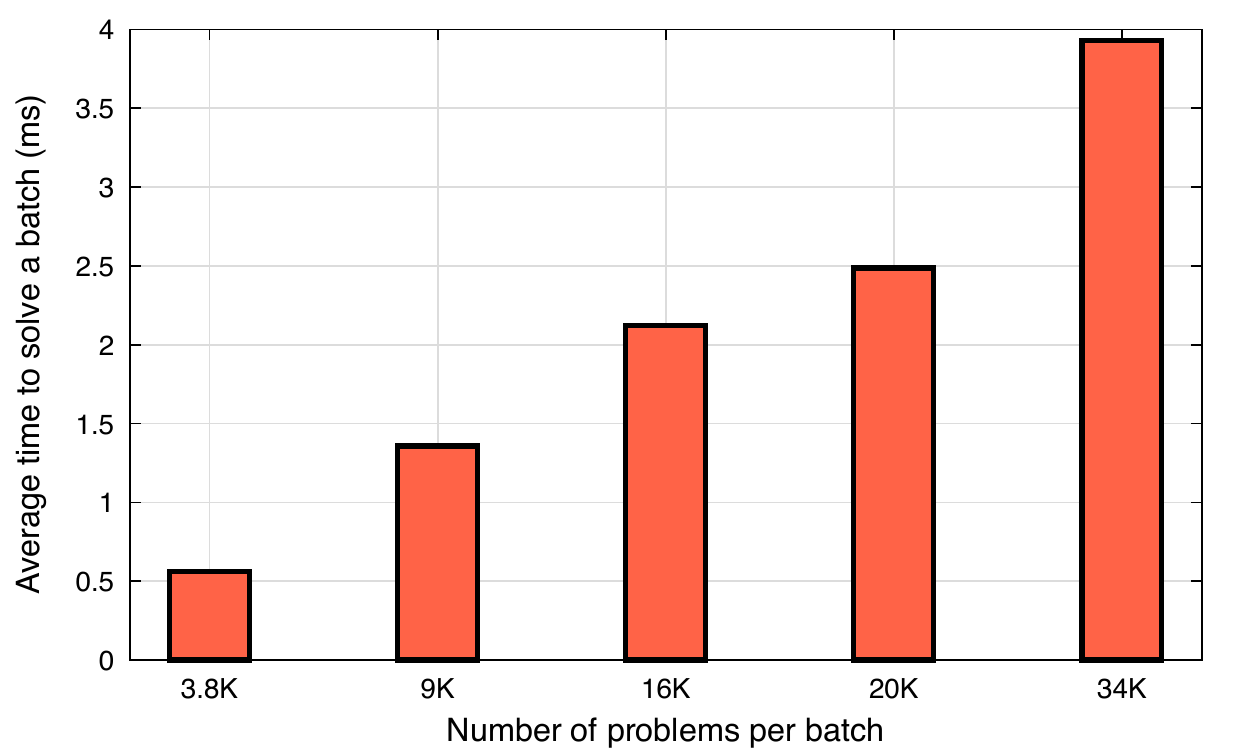}
    \caption{Performance of \ExaTron{} on a single GPU}
    \label{fig:single-gpu}
\end{figure}

\Cref{fig:single-gpu} depicts the average solution time of \ExaTron{} for different sizes of batches of branch subproblems listed in~\cref{tbl:exp-data}.
The time on the y-axis is the average computation time in milliseconds taken by \ExaTron{} to solve each batch within an ADMM iteration.

As illustrated in the figure, the performance of \ExaTron{} generally scales linearly with respect to the batch size.
This is expected since \ExaTron{} solves the subproblems inside a batch in parallel, meaning that increasing the batch size would linearly increase its computation time as well.
Moreover, all the subproblems share the same formulation~\cref{eq:exp-branch},  differing only by the parameter values.

We note that there are two main factors contributing to the computation time of a batch:
(i) the number of subproblems per batch (i.e., batch size) and
(ii) the average time for \ExaTron{} to solve each subproblem in a batch.
The latter is related to the level of difficulty of the subproblems in a batch.
Although two batches are of the same size, one batch can show more computation time than the other if the subproblems in that batch are more difficult to solve than those in the other.
In our case, the difficulty of batches was not significantly different from each other so \ExaTron{} showed a linear scaling over them.
However, this is not always the case, as demonstrated in~\cref{subsec:exp-multiple-gpu}, where we present a load imbalance on multiple GPUs that was caused by different level of difficulties among batches.

\subsection{Performance on multiple GPUs}
\label{subsec:exp-multiple-gpu}

\begin{figure}[htbp!]
    \centering
    \includegraphics[scale=.65]{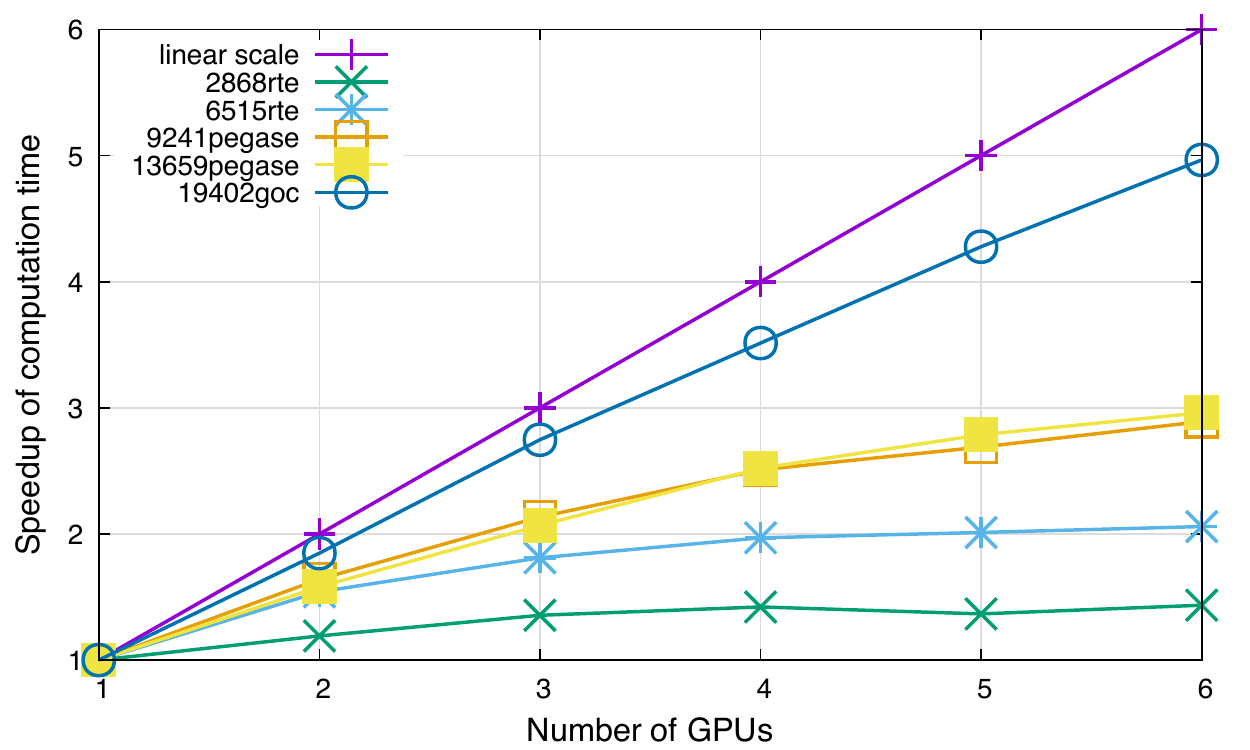}
    \caption{Performance of \ExaTron{} on multiple GPUs}
    \label{fig:multiple-gpu}
\end{figure}


We can achieve greater speedups of \ExaTron{} by employing multiple GPUs.
The GPU-aware message passing interface (MPI) with Nvidia's GPUDirect supports direct GPU-GPU communication by allowing GPUs to directly send or receive data from and to the device memories, without staging through host memory.
We employ such a GPU-aware MPI to communicate between GPUs that are connected via NVLink on Summit.
We note that the MPI communication can become more expensive when we start to use more than 3 GPUs, since this involves cross-socket communications.
This degradation of communication time will be observed in the experiments later in this section.

\Cref{fig:multiple-gpu} shows the speedup of \ExaTron{} when we parallelize the computation across different GPUs (up to the 6 GPUs available on a node in the Summit supercomputer).
Branch problems are evenly dispatched among 6 MPI processes in the order of branch indices, and the speedup is computed based on the timing of the root process.
At each ADMM iteration, the root process distributes variable values over multiple GPUs, solves its own batch, and gathers variable values back from other GPUs.
Therefore, the timing of the root process represents the synchronous time for all the MPI processes to finish their own solves, and communicate their solutions back to the root process.

As expected, we obtained a larger speedup as we increase the number of GPUs.
In the case of 19402goc, it shows almost a perfect linear scaling:
we obtained 5 times faster computation time when we used 6 GPUs.
In general, a larger speedup is achieved for larger test instances, because the speedup is related to the maximum number, say $N$, of subproblems that a single GPU can solve in parallel.
Therefore, the size of a batch divided by $N$ will determine the upper bound on the number of GPUs that we can benefit from employing them.
This explains why the speedup is saturated for smaller batches, for example 2868rte and 6515rte.
We also note that the slope of the changes of the speedup was slightly decreased when we used more than 3 GPUs.
We think this is because of the increased communication cost for cross-socket communications.

In addition to the size of a batch, another factor contributing to the speedup is the load balance among GPUs.
As we briefly discussed in~\cref{subsec:exp-single-gpu}, the level of difficulty of problems affect the computation time of a batch.
Although each GPU is assigned to a batch of the same size, the speedup may be degraded if the computational load is imbalanced in terms of the level of difficulty between GPUs.
In our case, the level of difficulty of a batch is measured by its solution time on a GPU.
Since the root process operates synchronously, some GPU may finish its computation much earlier than others, making it idle until all the other GPUs finish their work.

We quantify this load imbalance of our data using the metrics described 
in~\cite{Pearce12,kim2019asynchronous} and present their values in~\cref{tbl:load-imbalance}.
The percent imbalance of a problem instance $p$---where $p$ is divided into batches on GPUs in our case---at iteration $k$ is defined as
\begin{equation}
    v_{pk} := \left(\frac{t^{\max}_{pk}}{\bar{t}_{pk}} - 1 \right) \times 100\%,
\end{equation}
where $t^{\max}_{pk}$ and $\bar{t}_{pk}$ are the maximum and mean computation times of a batch among GPUs at iteration $k$.
We define $\overline{\nu}_p = \max_k \nu_{pk}$, $\underline{\nu}_p=\min_k\nu_{pk}$, and $\nu^m_p := \text{mean}_k \nu_{pk}$.
Hence, a smaller $\nu^m_p$ implies that the load is balanced better than the case with a larger $\nu^m_p$.

\Cref{tbl:load-imbalance} clearly shows that the workload of 19402goc among GPUs is much better balanced than the others, providing another insight into its superior speedup.
We note that the load balance of 13659pegase is worse than that of 9241pegase.
This explains why 13659pegase shows a very small increase of the speedup compared to 9241pegase, although its size is about 25\% larger than that of 9241pegase.

\begin{table}[htbp!]
    \centering
    \caption{Load imbalance metric values}
    \label{tbl:load-imbalance}
    \begin{tabular}{|r|r|r|r|}
    \hline
        \multicolumn{1}{|c|}{Data ($p$)} & $\overline{\nu}_p$ & $\underline{\nu}_p$ & $\nu^m_p$\\\hline
        2868rte & 481.58\% & 2.85\% & 46.76\%\\
        6515rte & 480.88\% & 0.54\% & 37.02\%\\
        9241pegase & 475.48\% & 1.40\% & 32.79\%\\
        13659pegase & 469.64\% & 5.98\% & 45.76\%\\
        19402goc & 469.79\% & 2.14\% & 9.04\%\\
        \hline
    \end{tabular}
\end{table}

\begin{figure}[htbp!]
    \centering
    \includegraphics[scale=.4]{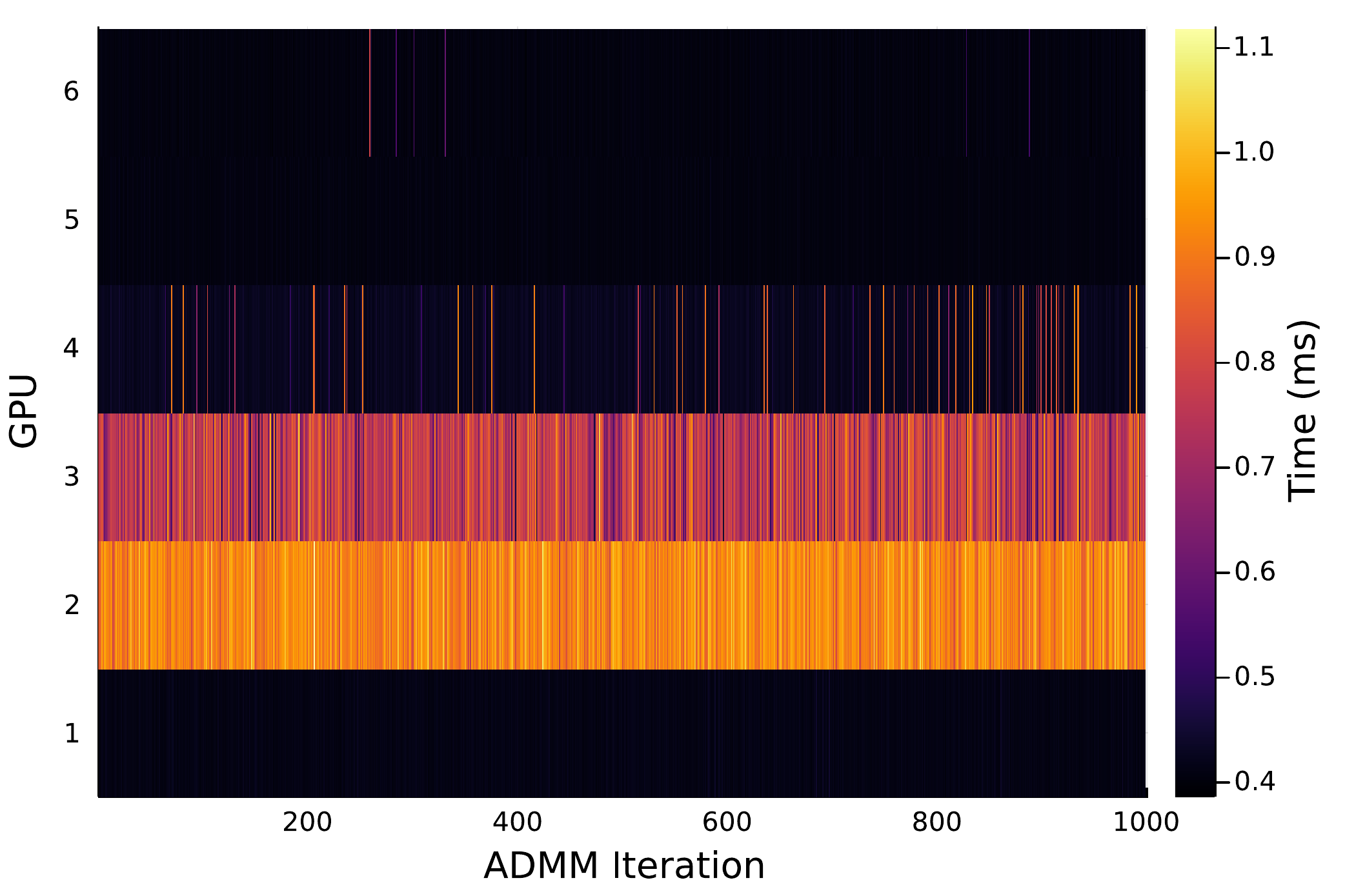}
    \caption{Load imbalance among GPUs for 13659pegase}
    \label{fig:load-imbalance}
\end{figure}

\Cref{fig:load-imbalance} presents a verification of the load imbalance of 13659pegase problem.
It depicts a heatmap where we sample 1,000 ADMM iterations for visibility and the value corresponds to the computation time in milliseconds of a batch that is assigned to each GPU.
As we observe in the figure, the computation time of the second and the third GPU was almost twice more than the others throughout the iterations.
An asynchronous solve combined with a load balancing scheme is our future research topic to alleviate this load imbalance issue.


\subsection{Performance comparison: CPU vs GPU}
\label{subsec:exp-cpu}

\begin{figure}[htbp!]
    \centering
    \includegraphics[scale=.65]{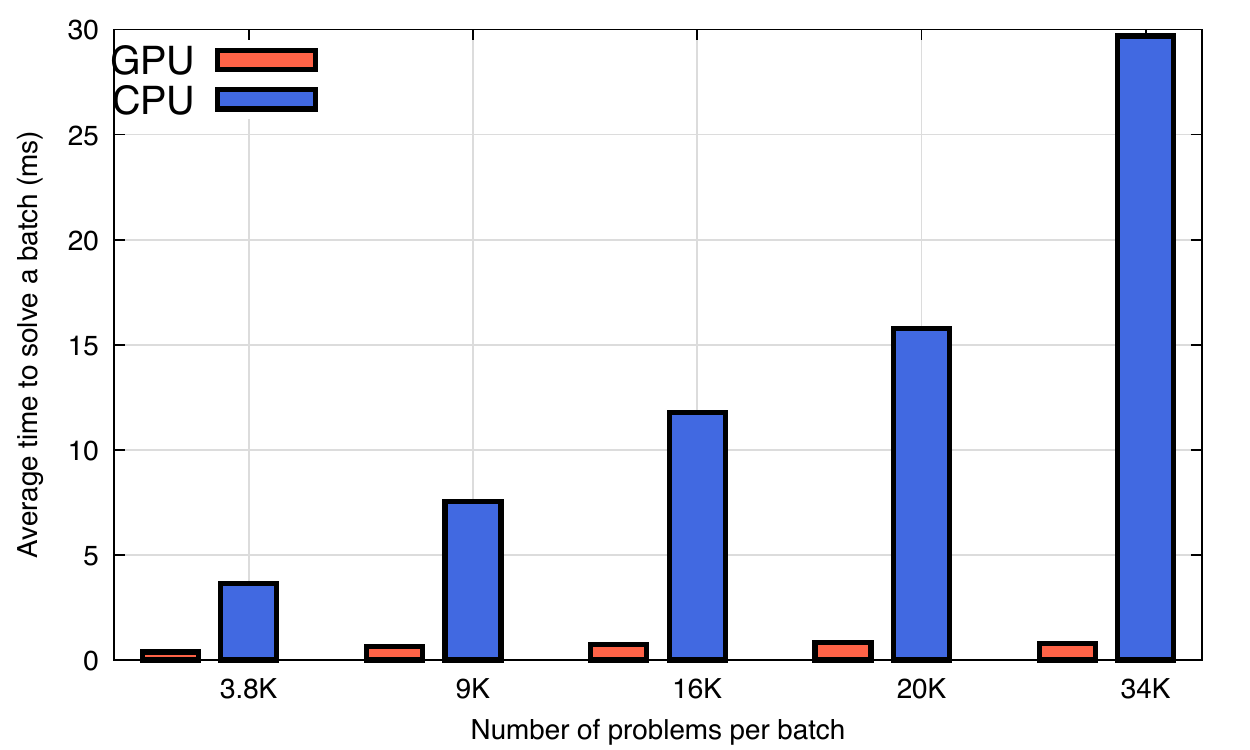}
    \caption{Performance comparison: 40 CPUs vs 6 GPUs on a single Summit node}.
    \label{fig:exp-cpu-vs-gpu}
\end{figure}

We next compare the performance of \ExaTron{} between a parallel CPU implementation and the GPU implementation.
This experiment was run on a single Summit node with 6 GPUs and 40 CPUs.
For the CPU run, we use the MPI library to implement the parallel communication between the CPU processes.
Similar to the experiments for multiple GPUs in~\cref{subsec:exp-multiple-gpu}, we measured the timing of the root process for CPUs that includes the cost for synchronizing distributed solves over the 40 cores.

Since \ExaTron{} can operate in either CPU or GPU mode, the same \ExaTron{} package was used to test the CPU implementation.
We note that both modes implement the same algorithm described in Algorithm~\cref{alg:exatron}, hence the sequence of calling their functions is identical.
The only differences are in the implementation of such functions---especially functions listed in~\cref{fig:core-ops}---where we follow our kernel design principles described in~\cref{sec:impl} for the GPU implementation.

In~\cref{fig:exp-cpu-vs-gpu}, the computation time of the CPU implementation shows a linear increase of
with respect to the batch size.
However, the average computation time increases faster than that of the GPU implementation: the computation time of \ExaTron{} on 6 GPUs is about 9--35 times faster than the CPU implementation using 40 cores.
Most of the speedup relates to the GPU's massive parallel computation capability.

\section{Conclusion}
\label{sec:conclusion}

Large-scale nonlinear programming can be tackled through Lagrangian
decomposition.
Such decomposition results in a batch of nonlinear programming problems
to solve at each iteration of its algorithm.
We have developed \ExaTron{} for efficient batched nonlinear programming
using GPUs.
It implements a trust-region Newton algorithm for bound constrained nonlinear
programming problems and works fully on GPUs without data transfers between
CPU and GPU.
This removes expensive data transfer cost which could be significant
especially when the size of problems in the batch is small.
We presented our design principles and implementation details for our
kernel function for efficient utilization of GPUs.
Experimental results over large-scale ACOPF problems decomposed into
components through ADMM algorithm provided linear scaling of computation
speed of \ExaTron{} with respect to the batch size and the number of GPUs.
On a single Summit node, the algorithm running on GPUs achieved more than 35 times speedup than on CPUs.

We conclude this paper by discussing several directions of future work. First of all, we plan to apply the ADMM algorithm with our GPU batch solver for solving multi-period multi-scenario optimal power flow problem. We found that only minor modifications are required to extend the current ADMM algorithm for solving the problem, while requiring more computing resources (i.e., multiple nodes and GPUs). Moreover, as we introduce more GPUs, a better design for MPI communication will be required. In particular, we have already observed that the computational load can be significantly imbalanced over multiple processes. Advanced asynchronous algorithms (e.g.,~\cite{kim2019asynchronous}) will be required to alleviate the load imbalance.


\bibliographystyle{siamplain}
\bibliography{references}

\vspace{0.3in}
\noindent\fbox{\parbox{0.95\textwidth}{
The submitted manuscript has been created by UChicago Argonne, LLC, Operator of Argonne National Laboratory (``Argonne''). Argonne, a U.S. Department of Energy Office of Science laboratory, is operated under Contract No. DE-AC02-06CH11357. The U.S. Government retains for itself, and others acting on its behalf, a paid-up nonexclusive, irrevocable worldwide license in said article to reproduce, prepare derivative works, distribute copies to the public, and perform publicly and display publicly, by or on behalf of the Government. The Department of Energy will provide public access to these results of federally sponsored research in accordance with the DOE Public Access Plan (http://energy.gov/downloads/doe-public-access-plan).}
}

\end{document}